\documentclass [12pt]  {article} 
\usepackage{amssymb}
\def \ZZ{{\mathbb{Z}}}
\def \QQ{{\mathbb{Q}}}

\def \CC{{\mathbb{C}}}
\def \FF{{\mathbb{F}}}

\begin{document}

\begin{center}
{\Large {\bf Some remarks on bielliptic and trigonal curves}}\\
\bigskip
{\tiny {\rm BY}}\\
\bigskip
{\sc Andreas Schweizer\footnote{The author was supported by 
the National Research Foundation of Korea (NRF) grant funded
by the Korean government (MSIP) (ASARC, NRF-2007-0056093).}}\\
\bigskip
{\small {\rm Department of Mathematics,\\
Korea Advanced Institute of Science and Technology (KAIST),\\ 
Daejeon 305-701\\
South Korea\\
e-mail: schweizer@kaist.ac.kr}}
\end{center}
\begin{abstract}
\noindent
We prove some results on algebraic curves $X$ of genus $g\geq 2$ in 
characteristic $0$. For example: Assume that $X$ has an automorphism
$\sigma$ of prime order $p\geq 5$. If $\sigma$ has no fixed points, 
then $X$ cannot be trigonal. On the other hand, if $\sigma$ has fixed 
points, then $X$ is bielliptic only if it belongs to one of three 
extremal types of curves of small genus.   
\end{abstract}

\subsection*{1. Introduction and basic facts}

Unless said otherwise, in this paper curve means a connected, smooth,
projective, algebraic curve over the complex numbers. Equivalently, 
one can think of a compact Riemann surface. 
\par
In [M], using the uniformization of compact Riemann surfaces of genus 
$g\geq 2$ by Fuchsian groups, Maclachlan proved, among other results, 
the following theorem.
\\ \\
{\bf Theorem 1.1. (Maclachlan)} [M, Theorem 2] \it 
Let $X$ be a hyperelliptic curve and $G$ a subgroup of $Aut(X)$ 
such that the covering $X\to X/G$ is totally unramified. 
Then $|G|$ divides $4$ and $G$ has exponent $2$.
\rm
\\ \\
We first give a short and elementary proof of this theorem. 
Then we work out similar and related results for other classes
of curves, notably for bielliptic and trigonal curves.
\\ \\
{\bf Proof.} \rm 
Denote the genus of $X$ by $g$. The hyperelliptic involution $\tau$ 
lies in the center of $Aut(X)$. So $G$ acts on the $2g+2$ fixed points 
of $\tau$, and as the action of $G$ is without fixed points, $|G|$ 
divides $2g+2$. On the other hand, since the covering is unramified, 
the Hurwitz formula implies that $|G|$ divides $2g-2$. Subtraction 
already shows that $|G|$ divides $4$.
\par
With an eye towards generalization, we prove the second claim 
independently of the first by counting the ramified points on 
$X$ of the covering $X\to X/\langle\tau,G\rangle$.
\par
The subcover $X/\langle G\rangle\to X/\langle\tau,G\rangle$ is the
hyperelliptic one and has $2h+2$ ramified points, where $h$ is the
genus of $X/\langle G\rangle$. They all split completely in $X$, 
which gives $|G|(2h+2)=2g+2+4|G|-4$ ramified points, each with 
ramification index $2$. 
\par
Taking the $2g-2$ ramified points of $X\to X/\langle\tau\rangle$ 
into account, there must still be $2|G|-2$ ramified points in 
$X/\langle\tau\rangle\to X/\langle\tau,G\rangle$, each with 
ramification index $2$ and splitting into $2$ points in $X$. 
Note that the covering 
$X/\langle\tau\rangle\to X/\langle\tau,G\rangle$ is Galois with 
group $\widetilde{G}$ induced by (and isomorphic to) $G$, but not 
necessarily fixed-point-free. The stabilizer in $\widetilde{G}$ 
of each of these $2|G|-2$ points on $X/\langle\tau\rangle$ is an 
involution. Since $X/\langle\tau\rangle$ has genus $0$, each such 
involution has exactly $2$ fixed points. So $\widetilde{G}$ 
(and hence $G$) must have at least $|G|-1$ involutions, 
i.e. every non-trivial element is an involution.
\hfill$\Box$
\\ \\
Maclachlan's theorem says that hyperelliptic curves cannot have 
automorphisms of certain types, for example fixed-point-free of 
prime order $p\geq 3$. In the sequel we want to generalize this 
to other types of curves. 
\par
On the other hand, it is clear that the curve $Y=X/G$ in Theorem 
1.1 must also be hyperelliptic. So another interpretation is that 
certain coverings $X$ of a hyperelliptic curve $Y$ cannot be
hyperelliptic. In this form we also want to generalize it, keeping
the curve $Y$ hyperelliptic, but taking other types of curves $X$.
\par
One generalization of a hyperelliptic involution is the notion of
a $b$-hyperelliptic involution, that is, an involution such that
the genus of the quotient curve is $b$. We concentrate on the 
case $b=1$; then the involution and the curve are called
{\bf bielliptic}. Some of the literature also uses the expression
elliptic-hyperelliptic for this.
\par
Another generalization of $X$ being hyperelliptic is the existence 
of an $n$-gonal map. This is a surjective morphism of degree $n$ 
from $X$ to a projective line. Note that for $n>2$ such a morphism 
is not necessarily Galois, i.e. does not necessarily come from an
automorphism of $X$. This makes it much more difficult to establish 
or exclude the existence of such a map. We are mainly interested in
the case $n=3$; then the map and the curve are called {\bf trigonal}.
\par
Throughout the paper we will freely use the following easy
consequences of the Castelnuovo inequality:
\begin{itemize}
\item The hyperelliptic involution is unique.
\item On a bielliptic curve of genus $g\geq 6$ the bielliptic 
involution is unique.
\item A curve of genus $g\geq  5$ has at most one trigonal map.
\item A hyperelliptic curve of genus $g\geq 3$ cannot be trigonal.
\item A hyperelliptic curve of genus $g\geq 4$ cannot be bielliptic.
\end{itemize}

\noindent
Finally, we'll need the following results.
\\ \\
{\bf Theorem 1.2.} \it 
Let $Y$ be a hyperelliptic curve of genus $h$, and let $X\to Y$ be
a smooth, cyclic covering of degree $n$.
\begin{itemize}
\item[(a)] {\bf (Bujalance)} {\rm [B, Theorem]}
If $n=2$, then $X$ is $b$-hyperelliptic for some $b$ with
with $b\leq\lfloor\frac{h-1}{2}\rfloor$. 
\par 
In particular, a smooth degree $2$ cover of a genus $2$ curve 
is a genus $3$ curve that is hyperelliptic (and bielliptic).
\item[(b)] {\bf (Accola)} {\rm [A, Lemma 2]} 
If $n$ is odd, then the hyperelliptic 
involution of $Y$ lifts to $n$ different involutions on $X$ that 
are all $b$-hyperelliptic where $b=\frac{(n-1)(h-1)}{2}$. Together
with the automorphism of order $n$ from $Gal(X/Y)$ they generate
a dihedral group $D_n$ of order $2n$.
\par
In particular, a smooth degree $3$ Galois cover of a genus $2$ curve 
is a genus $4$ curve that has at least $3$  bielliptic involutions.
\end{itemize}
\rm

\subsection*{2. Bielliptic curves}

We start with a quick corollary to the known results from the 
previous section.
\\ \\
{\bf Corollary 2.1.} \it
Let $X$ be a curve of genus $5$. If $X$ has an automorphism 
$\sigma$ of order $4$ such that $\sigma$ has no fixed points, 
then $X$ is bielliptic.
\rm
\\ \\
{\bf Proof.} \rm 
As $X/\langle\sigma^2\rangle$ is an unramified degree $2$ cover 
of $X/\langle\sigma\rangle$, the genus of $X/\langle\sigma^2\rangle$
can only be $1$ (if $\sigma^2$ has fixed points) or $3$. In the
first case we are done. 
\par
In the second case we get a chain 
$X\to X/\langle\sigma^2\rangle\to X/\langle\sigma\rangle$
of smooth coverings of degree $2$ with genera $5$, $3$, and 
$2$. By Theorem 1.2 (a) the genus $3$ curve is hyperelliptic, 
and the genus $5$ curve is hyperelliptic or bielliptic. But 
it cannot be hyperelliptic by Theorem 1.1.
\hfill$\Box$
\\ \\
Ultimately Theorem 1.2 and Corollary 2.1 are based on lifting the 
hyperelliptic involution to an unramified cover. Our next result 
requires lifting the hyperelliptic involution to a ramified cover, 
which is a much more tricky problem. The paper [CT] discusses
conditions under which this is possible. Compared to our previous
statements we change the names of the curves and their automorphisms 
so that they fit precisely with those in the somewhat technical 
conditions in [CT].
\\ \\
{\bf Proposition 2.2.} \it
Let $Y$ be a non-hyperelliptic curve of genus $4$. If $Y$ has 
an automorphism $\tau$ of order $4$, then $Y$ is bielliptic.
\rm
\\ \\
{\bf Proof.} \rm 
Since $Y$ is non-hyperelliptic, the curve $X=Y/\langle\tau^2\rangle$
can only have genus $1$ or $2$. In the first case we are done. 
\par
In the second case $\tau^2$ has exactly $2$ fixed points $P$ and $Q$ 
on $Y$. From the Hurwitz formula we see that $Y/\langle\tau\rangle$ 
has genus $1$ and $\tau$ also fixes $P$ and $Q$. So $\widetilde{\tau}$,
the involution induced by $\tau$ on $X$ is different from $\sigma$,
the hyperelliptic involution of $X$. 
\par
As $\sigma$ commutes with $\widetilde{\tau}$, it acts on its fixed 
points $\widetilde{P}$ and $\widetilde{Q}$, the images of $P$ and $Q$
on $X$. Moreover, different involutions must have disjoint fixed points.
Hence $\sigma(\widetilde{P})=\widetilde{Q}$. Since the Weierstrass
points of $X$ are exactly the fixed points of $\sigma$, we see that 
$\widetilde{P}$ and $\widetilde{Q}=\sigma(\widetilde{P})$, the only two 
points on $X$ that ramify in $Y$, are not Weierstrass points of $X$. So 
our constellation satifies the conditions of [CT, Theorem 2.1], which 
allows us to conclude that $\sigma$ lifts to an automorphism of $Y$.
\par
Actually, $\sigma$ cannot lift to an automorphism $\psi$ of order 
$4$, because then we would have $\psi^2=\tau^2$ and as before 
$X/\langle\sigma\rangle=Y/\langle\psi\rangle$ would have genus $1$. 
So $\sigma$ lifts to an involution $\psi$ on $Y$ that commutes with 
$\tau^2$. Then the other two intermediate curves between $Y$ and the 
genus zero curve $Y/\langle\tau^2,\psi\rangle$ have genera $2$ and 
$0$ (which would contradict non-hyperellipticity of $Y$) or $1$ and 
$1$, which implies that in this case $Y$ actually has at least $2$ 
bielliptic involutions.
\hfill$\Box$
\\ \\
Here is an analogue of Maclachlan's result for bielliptic curves.
\\ \\
{\bf Proposition 2.3.} \it
Let $Y$ be a hyperelliptic curve of genus at least $4$ and let
$X\to Y$ be a smooth Galois cover with Galois group $G$ such that
$X$ is bielliptic. Then
\begin{itemize}
\item[(a)] At least half of the elements of $G$ are involutions.
\item[(b)] If $G$ is abelian, then $G$ has exponent $2$ and $|G|$ 
divides $8$.
\end{itemize}
\rm

\noindent
{\bf Proof.} \rm 
\par
(a) Since $g(X)\geq 6$, the bielliptic involution $\tau$ is unique
and hence commutes with $G$. As $Y$ cannot be bielliptic, $\tau$
induces the hyperelliptic involution on $Y$. Counting ramified points
in the covering $X\to X/\langle\tau,G\rangle$ as in the proof of the
second claim of Theorem 1.1, one sees that $G$ must have at least
$|G|/2$ involutions.
\par
(b) If $G$ is abelian, part (a) implies that $G$ must be $2$-elementary 
abelian. In general, $G$ is isomorphic to a subgroup of 
$Aut(X/\langle\tau\rangle)$. But the biggest $2$-elementary abelian 
group contained in the automorphism group of a genus one curve is 
$\ZZ/2\ZZ\oplus\ZZ/2\ZZ\oplus\ZZ/2\ZZ$.
\hfill$\Box$
\\ \\
If in Proposition 2.3 we put the condition that the curve $Y$ is 
bielliptic instead of hyperelliptic, then the number of fixed 
points of the bielliptic involution fits with the Hurwitz formula 
for smooth covers, and the smoothness of the cover fits with the fact
that a surjective map between two curves of genus $1$ is unramified.
In short, the method of proof does not seem to imply any restrictions.
\par
On the other hand, for hyperelliptic $Y$ we can even say something
about the non-smooth case. But we need a little preparation.
\\ \\
{\bf Lemma 2.4.} \it
Let $X$ be a bielliptic curve of genus $g\geq 6$ with an automorphism 
$\sigma$ of prime order $p\geq 3$.
\begin{itemize}
\item[(a)] If $\sigma$ has fixed points, then necessarily $p=3$.
\item[(b)] If $\sigma$ has no fixed points, then the curve 
$X/\langle\sigma\rangle$ is also bielliptic.
\end{itemize}
\rm

\noindent
{\bf Proof.} \rm 
As $g\geq 6$, the bielliptic involution $\tau$ is unique. 
So $\tau$ and $\sigma$ commute. 
\par
For (a) we use that $\sigma$ induces an automorphism 
$\widetilde{\sigma}$ of order $p$ with fixed points on
the genus one curve $X/\langle\tau\rangle$, which is only
possible for $p=3$.
\par
For (b) note that any fixed point of $\widetilde{\sigma}$ on 
$X/\langle\tau\rangle$ would lift to one or two fixed points 
of $\sigma$ on $X$. So $\widetilde{\sigma}$ is fixed-point-free
and hence $X/\langle\tau, \sigma\rangle$ has genus $1$. This 
means that $\tau$ induces a bielliptic involution on 
$X/\langle\sigma\rangle$.
\hfill$\Box$
\\ \\
{\bf Proposition 2.5.} \it
Let $Y$ be a hyperelliptic curve of genus at least $4$. 
If $X\to Y$ is a Galois cover of degree relatively prime 
to $6$, then $X$ cannot be bielliptic.
\rm
\\ \\
{\bf Proof.} \rm 
Assume that $X$ is bielliptic. As $g(X)\geq 6$, the bielliptic 
involution $\tau$ is unique and lies in the center of $Aut(X)$. 
Moreover, because of $g(Y)\geq 4$ the Galois group $G$ of the 
covering $X\to Y$ cannot contain $\tau$. So $G$ induces an 
isomorphic group of automorphisms on $X/\langle\tau\rangle$.
As finite subgroups of automorphisms on a genus $1$ curve are 
solvable, $G$ is solvable. (Alternatively, we could argue here 
with the highly non-trivial fact that groups of odd order are 
solvable.) So there exists a chain of coverings 
$$X=X_1 \to X_2 \to \cdots \to X_s = Y$$
such that each covering $X_i \to X_{i+1}$ is Galois of prime degree.
By Lemma 2.4 (a) the biellipticity of $X_i$ implies that the 
covering $X_i \to X_{i+1}$ is unramified, and by Lemma 2.4 (b) this
implies that $X_{i+1}$ is bielliptic, too. By iteration we get that 
$Y$ is bielliptic, which cannot be for a hyperelliptic curve of 
genus bigger than $3$.
\hfill$\Box$
\\ \\
{\bf Remark 2.6.} 
With practically the same proof Proposition 2.5 holds under slightly more 
general conditions, for example if the degree of the Galois cover is 
odd and none of the ramification indices is divisible by $3$.
\\ \\
A bielliptic curve $X$ of genus $g(X)\leq 5$ can have more than one 
bielliptic involution. We recall some curves that reach the maximally
possible number of bielliptic involutions for the given genus.
\\ \\
The maximal possible number of automorphisms on a curve of genus $3$
is $168$. The unique genus $3$ curve that realizes this bound is the
modular curve $X(7)$. Its automorphism group is the simple group
$PSL_2(\FF_7)$. This group is also isomorphic to $GL_3(\FF_2)$. 
The projective model of this curve is 
$$x^3 y+y^3 z+z^3 x=0,$$
which is why it is also called the {\bf Klein quartic}.
Since its automorphism group is simple, this curve cannot be 
hyperelliptic. So by Theorem 1.2 (b) every involution must be
bielliptic. Thus it has $21$ bielliptic involutions, the maximum
possible for a curve of genus $3$.
\\ \\
By [CD, Corollary 6.9] a curve of genus $4$ can have at most $10$ 
bielliptic involutions, and there is exactly one curve of genus $4$ 
with $10$ bielliptic involutions. It is called {\bf Bring's curve}, 
and it is isomorphic to the modular curve $X_1(5,10)$. Its automorphism 
group is isomorphic to the symmetric group $S_5$. Actually, $120$ is 
also the maximum number of automorphisms for a curve of genus $4$.
\\ \\
By [KMV] a curve of genus $5$ can have $0$, $1$, $2$, $3$ or $5$ 
bielliptic involutions, and the genus $5$ curves with $5$ biellliptic 
involutions form a $2$-dimensional family. They are called 
{\bf Humbert curves}. Their automorphism groups have order $160$. This 
does not quite reach the maximally possible number of automorphisms 
on a curve of genus $5$, which is $192$.
\\ \\
For the proof of the next result we have to introduce yet another type 
of curve with an obvious automorphism of order $p$, namely for each 
prime $p\geq 5$ the {\bf Lefschetz curve}
$$y^p=x(x-1)$$
of genus $\frac{p-1}{2}$. Its full automorphism group is isomorphic to
$\ZZ/2p\ZZ$. The unique involution $x\mapsto 1-x$ is hyperelliptic,
and hence the curve is {\it not} bielliptic.
\\ \\
In [JKS2] we showed that there is only one bielliptic curve of genus 
$4$ with an automorphism of order $5$, namely Bring's curve. Now we
want to generalize this.
\\ \\
{\bf Theorem 2.7.} \it
Let $X$ be a bielliptic curve. Assume that $X$ has an automorphism 
$\sigma$ of prime order $p\geq 5$ with fixed points. Then one of the
following holds:
\begin{itemize}
\item[(a)] $g=5$, $p=5$ and $X$ is a Humbert curve;
\item[(b)] $g=4$, $p=5$ and $X$ is Bring's curve;
\item[(c)] $g=3$, $p=7$ and $X$ is the Klein quartic.
\end{itemize}
\rm

\noindent
{\bf Proof.} \rm 
By Lemma 2.4 (a) we must have $g\leq 5$. Moreover, $\sigma$ acts 
by conjugation on the bielliptic involutions. If it commutes with 
a bielliptic involution, it induces an automorphism of order $p$ 
with fixed points on the elliptic quotient curve, which is impossible.
Hence the number of biellitic involutions must be divisible by $p$.
\par
If $g=5$, the maximal number of bielliptic involutions is $5$ [KMV],
so necessarily $p=$ and $X$ has $5$ bielliptic involutions, i.e. $X$
is a Humbert curve.
\par
If $g=4$, the possible numbers of bielliptic involutions are 
$0$, $1$, $2$, $3$, $4$, $6$, $10$ [CD, p.600]. So $p=5$ 
and $10$ bielliptic involutions, i.e. $X$ is Bring's curve.
\par
If $g=3$, the possible prime divisors of $|Aut(X)|$ are $2$, $3$ and 
$7$. By [RR, Theorem 1] there are exactly two curves of genus $3$ with 
an automorphism of order $7$. One is the Klein quartic, the other one 
is the Lefschetz curve (which is not bielliptic).
\par
If $g=2$, the biggest possible prime divisor of $|Aut(X)|$ is 
$5$. By [RR, p.199] the only genus $2$ curve with an automorphism 
of order $5$ is the Lefschetz curve, which is not bielliptic.
\hfill$\Box$
\\

\subsection*{3. Trigonal curves}

We give another generalization of Maclachlan's theorem.
\\ \\
{\bf Theorem 3.1.} \it
Let $p$ be an odd prime. Let $X$ be a $p$-gonal curve of genus
$g>(p-1)^2$ and $G\leq Aut(X)$ a subgroup of fixed-point-free 
automorphisms. Then $|G|$ divides $p$.
\rm
\\ \\
{\bf Proof.} \rm 
Since $g>(p-1)^2$, by the Castelnuovo inequality the $p$-gonal 
map $\pi:X\to Y$ is unique. So $G$ induces an isomorphic group 
of automorphisms on the genus zero curve $Y$. Then $G$ acts without 
fixed points on the ramification points $P$ of $\pi$ on $X$. 
Consequently, they come in batches of size $|G|$ with the same 
ramification index $e_P$. Hence $|G|$ divides $\sum (e_P -1)$. 
Thus by the Hurwitz formula $|G|$ divides $2g-2+2p$.
\par
On the other hand, as the covering $X\to X/G$ is smooth, $|G|$ 
divides $2g-2$. Subtracting this, we see that $|G|$ divides $2p$. 
\par
So what is left is showing that every involution $\sigma$ on 
$X$ must have fixed points. Let $P$ be one of the two fixed 
points of the induced involution $\widetilde{\sigma}$ on $Y$.
Let $P_1,\ldots,P_r$ be the points on $X$ lying above $P$ and
let $e_1,\ldots,e_r$ be their ramification indices. Then 
$\sigma$ acts on these points. If $\sigma$ fixes no $P_i$, 
then they come in pairs with the same ramification index.
So $p=\sum_{i=1}^r e_i$ would be even, a contradiction.
\hfill$\Box$
\\ \\
Theorem 3.1 applies in particular to trigonal curves 
of genus $g\geq 5$.
\\ \\
{\bf Proposition 3.2.} \it
Let $Y$ be a hyperelliptic curve of genus $h\geq 3$,
and let $X\to Y$ be a Galois cover. Then $X$ cannot be 
trigonal.
\rm
\\ \\
{\bf Proof.} \rm 
If $X$ is trigonal, then, because of $g(X)\geq 5$, the trigonal 
map $X\to Z$ is unique. So $G$, the Galois group of $X\to Y$, 
induces an isomorphic group of automorphisms on $Z$. This implies
the existence of a trigonal map $Y\to Z/G$. But by Castelnuovo
$Y$ cannot be trigonal.
\hfill$\Box$
\\ \\
{\bf Corollary 3.3.} \it
Let $X\to Y$ be a smooth Galois cover of degree $n$ where $Y$ is 
a hyperelliptic curve of genus $h$. Then $X$ is trigonal if and only
if $n=3$, $h=2$, and $X$ has genus $4$.
\rm
\\ \\
{\bf Proof.} \rm 
Assume that $X$ is trigonal. If $g(X)\geq 5$, then $n=3$ by 
Theorem 3.1 and $h=2$ by Proposition 3.2; but this contradicts 
$g(X)\geq 5$. So for trigonal $X$ we must have $g(X)\leq 4$. 
Now the Hurwitz formula leaves only the possibilities $g=4$, 
$n=3$, $h=2$ or $g=3$, $n=2$, $h=2$. But by Theorem 1.2 (a)
the second possibility would imply the contradiction that 
$X$ is hyperelliptic.
\par
Conversely, a curve of genus $4$ is either hyperelliptic or 
trigonal. But by Theorem 1.1 a smooth Galois cover of degree 
$3$ cannot be hyperelliptic.
\hfill$\Box$
\\ \\
Concerning involutions on trigonal curves we have the following result.
\\ \\
{\bf Lemma 3.4.} \it
Let $X$ be a trigonal curve of genus $g$ and $\sigma$ an involution
on $X$.
\begin{itemize}
\item[(a)] If $g$ is odd, then $\sigma$ has exactly $4$ fixed points.
\item[(b)] If $g$ is even, then $\sigma$ has $2$ or $6$ fixed points.
\end{itemize}
\rm

\noindent
{\bf Proof.} \rm 
If $g\geq 5$, the trigonal map $X\to Y$ is unique, and hence 
$\sigma$ induces an involution $\widetilde{\sigma}$ on $Y$. 
Since $\widetilde{\sigma}$ has exactly two fixed points on $Y$, 
$\sigma$ can have at most $6$ fixed points on $X$. By the Hurwitz
formula this means $2$ or $6$ fixed points if $g$ is even, and
$0$ or $4$ if $g$ is odd. But $0$ is excluded by Theorem 2.1.
\par
If $g=4$, then a priori there could be $2$, $6$ or $10$ fixed 
points. But $10$ fixed points would mean that the involution
is hyperelliptic, which is not possible for a trigonal curve
of genus bigger than $2$. 
\par
For the same reason $8$ fixed points are not possible if $g=3$.
But $0$ fixed points are not possible either for $g=3$, as by
Theorem 1.2 this would also imply that $X$ is hyperelliptic.
\hfill$\Box$
\\ \\
{\bf Corollary 3.5.} \it
Let $X$ be a curve of genus $g\equiv 1\ mod\ 4$. If $Aut(X)$ has 
a subgroup $H$ isomorphic to $\ZZ/2\ZZ\oplus\ZZ/2\ZZ$, then $X$ 
cannot be trigonal.
\rm
\\ \\
{\bf Proof.} \rm 
Assume that $X$ is trigonal. Then by Lemma 3.4 each of the $3$ 
involutions in $H$ has exactly $4$ fixed points. But then applying
the Hurwitz formula to the covering $X\to X/H$ would give 
$g\equiv 3\ mod\ 4$.
\hfill$\Box$
\\ \\
The famous Theorem of Faltings says that a curve $X$ of genus $g\geq 2$ 
over a number field $K$ has only finitely many $K$-rational points.
Likewise, it has only finitely many $L$-rational points for any fixed
cubic extension $L$ of $K$. 
But since $K$ has infinitely many cubic extensions, $X$ could have 
infinitely many cubic points in total. 
Here a cubic point over $K$ means a point $P$ that is $L$-rational 
for a cubic extension $L$ of $K$ where $L$ might depend on $P$.
\par
For example, if $X$ is trigonal over $K$, i.e. if there is a degree 
$3$ cover $X\to Y$, defined over $K$, such that $Y$ is a projective 
line over $K$, then $X$ has infinitely many cubic points over $K$, 
namely at least one over each of the infinitely many $K$-rational 
points of $Y$. 
\par
This is the reason why trigonality is one of the key points to 
examine, if for a given ensemble of curves one wants to decide 
which of them have infinitely many cubic points. See for example
[JKS1], where we investigated certain modular curves under this 
aspect to determine which finite groups occur infinitely often 
as torsion groups of elliptic curves over cubic number fields.
\par
Now a curve of genus $2$ is always trigonal. If $P$ is any 
non-Weierstrass point, there is a rational function with a triple
pole at $P$. But the question is whether there is a $K$-rational 
trigonal map. By what was just said, a $K$-rational point that is 
not a Weierstrass point would suffice for this. But not every 
genus $2$ curve has such a point. In [JKS1, Lemma 2.1] we also 
showed that if a genus $2$ curve has $3$ or more $K$-rational
points (Weierstrass or not), then it is trigonal over $K$, and
hence has infinitely many cubic points over $K$.
\par
Considering points that are genuinely cubic over $K$ (i.e. cubic
over $K$ but not $K$-rational), we now prove that a curve of genus 
$2$ either has infinitely many such points or none. 
\\ \\
{\bf Theorem 3.6.} \it
For a curve $X$ of genus $2$ over a number field $K$ the following
are equivalent:
\begin{itemize}
\item[(a)] $X$ has a cubic point $P$ over $K$ that is not $K$-rational.
\item[(b)] $X$ has a trigonal map that is defined over $K$.
\item[(c)] $X$ has infinitely many cubic points over $K$.
\end{itemize}
\rm

\noindent
{\bf Proof.} \rm 
Let $L$ be the cubic extension of $K$ generated by $P$. Let $\iota_1$,
$\iota_2$, $\iota_3$ be the three different $K$-embeddings of $L$ into 
$\CC$. 
\par
If there is a $\overline{K}$-rational function on $X$ with pole
divisor $\iota_i (P)+\iota_j (P)$ ($i\neq j$), this divisor must lie 
in the canonical class. We can apply a suitable automorphism from 
$Gal(\overline{K}/K)$ that permutes the three embeddings cyclically, 
and get that $\iota_j (P)+\iota_k (P)$ also lies in the canonical class.
But then $\iota_i (P)-\iota_k (P)$ would be a principal divisor, which
is impossible on a curve of positive genus. 
\par
Thus the Riemann-Roch space of the $K$-rational divisor 
$D=\iota_1 (P)+\iota_2 (P)+\iota_3 (P)$ contains a $K$-rational function
whose pole divisor actually is $D$. This function is a $K$-rational
trigonal map.
\par
The conclusion $(b)\to (c)$ is clear, and $(c)\to (a)$ is trivial.
\hfill$\Box$
\\ \\
Finally we show that among the genus $2$ curves that have no 
$\QQ$-rational points some are trigonal over $\QQ$ and some 
are not.
\\ \\
{\bf Example 3.7.}
The genus $2$ curve
$$y^2 =-x^6 -1$$
has no $\QQ$-rational points and no cubic points over $\QQ$. In 
fact, it has no real points at all. The points lying over the point 
at infinity of the $x$-line are not real, as $-1$ is a square in their 
residue field. 
\par
Admittedly, the degree $3$ map from the curve to $C:\ y^2=-x^2 -1$ 
is defined over $\QQ$. But $C$ is a genus $0$ curve without rational 
points. So over $\QQ$ this map is not a trigonal map. Over $\QQ(i)$,
where $C$ is a projective line, the map will of course be trigonal.
\\ \\
{\bf Example 3.8.}
The genus $2$ curve
$$y^2 =-x^6 +13$$
also has no $\QQ$-rational points. Here we are using that the elliptic 
curve $y^2 =x^3 +13$ has rank $0$ and trivial torsion, and that the
points at infinity are not $\QQ$-rational for exactly the same reasons 
as in the previous example.
\par
Again, there is an obvious map of degree $3$, defined over $\QQ$, to 
the genus $0$ curve $y^2 =-x^2 +13$. But this curve has $\QQ$-rational 
points, for example $(\pm 2,\pm 3)$ or $(\pm 3,\pm 2)$, so it is a 
projective line over $\QQ$, and the degree $3$ map to it is a trigonal 
map over $\QQ$.
Indeed, $(\pm \sqrt[3]{2}, \pm 3)$ and $(\pm \sqrt[3]{3}, \pm 2)$ are
genuine cubic points on our genus $2$ curve.
\\

\subsection*{\hspace*{10.5em} References}
\begin{itemize}

\item[{[A]}] R.~Accola: On lifting the hyperelliptic involution,
\it Proc. Amer. Math. Soc. \bf 122 \rm (1994), 341-347

\item[{[B]}] E.~Bujalance: A classification of unramified double 
coverings of hyperelliptic Riemann surfaces, \it Arch. Math. (Basel) 
\bf 47 \rm (1986), 93-93

\item[{[CD]}] G.~Casnati and A.~Del Centina:
The moduli spaces of bielliptic curves of genus 4 with more bielliptic 
structures, \it J. London Math. Soc. (2) \bf 71 \rm (2005), 599-621 

\item[{[CT]}] A.~F.~Costa and P.~Turbek: Lifting involutions to 
ramified covers of Riemann surfaces, \it Arch. Math. (Basel) \bf 81 
\rm (2003), 161-168

\item[{[JKS1]}] D.~Jeon, C.~H. Kim, and A.~Schweizer: On the torsion of
elliptic curves over cubic number fields, \it Acta. Arith. \bf 113 \rm
(2004), 291-301

\item[{[JKS2]}] D.~Jeon, C.~H. Kim, and A.~Schweizer: Bielliptic 
intermediate modular curves, \it preprint\rm

\item[{[KMV]}] T.~Kato, K.~ Magaard, H.~V\"olklein:
Bi-elliptic Weierstrass points on curves of genus $5$,
\it Indag. Math. (N.S.) \bf 22 \rm (2011), 116-130 

\item[{[M]}] C.~Maclachlan: Smooth coverings of hyperelliptic surfaces,
\it Quart. J. Math. Oxford Ser. (2) \bf 22 \rm (1971), 117-123

\item[{[RR]}] G.~Riera and R.~Rodriguez: Riemann surfaces and abelian 
varieties with an automorphism of prime order, \it Duke Math. J. \bf 69 
\rm (1993), 199-217

\end{itemize}

\end{document}